\documentclass[a4paper,12pt]{article}

\usepackage{latexsym,amsmath,amssymb,theorem,color }
\usepackage{xypic,amsmath}
\xyoption{all}
\usepackage{geometry}
\usepackage{enumerate}

\newcommand{\C}{\mathsf{C}}
\newcommand{\Proof}{\noindent\textbf{Proof }}
\newcommand{\nsq}{\square \hspace{-2.9mm}\raisebox{0.75mm}{\scriptsize $\diagdown$ }}

\newtheorem{examples}{Example}[section]
{\theorembodyfont{\rmfamily}}
\newtheorem{prop}[examples]{Proposition}
\newtheorem{theorem}[examples]{Theorem}
{\theorembodyfont{\rmfamily}}
\newtheorem{Corollary}[examples]{Corollary}
\newtheorem{lemma}[examples]{Lemma}

\def\w{$\omega$}
\def\o{\circ}
\def\bpartial{\boldsymbol{\partial}}
\def\s{\mathsf s}

\def\bvarepsilon{\boldsymbol{\varepsilon}}
\def\bGamma{\boldsymbol{\Gamma}}
\newcommand{\xdirects}[2]{\quad \def\objectstyle{\scriptstyle} \objectmargin={0pt}
\xy
(0,0)*+{}="a",(0,-6)*+{\rule{0em}{1.5ex}#2}="b",(7,0)*+{\;#1}="c"
\ar@{->} "a";"b" \ar @{->}"a";"c"
\endxy }
\newcommand{\directs}[2]{\def\objectstyle{\scriptstyle} \objectmargin={0pt}
\xy
(0,4)*+{}="a",(0,-2)*+{\rule{0em}{1.5ex}#2}="b",(7,4)*+{\;#1}="c"
\ar@{->} "a";"b" \ar @{->}"a";"c" \endxy }

\newcommand{\sdirects}[2]{\def\objectstyle{\scriptstyle} \objectmargin={0pt}
\xy
(0,2.2)*+{}="a",(0,-2.5)*+{\rule{0em}{1.5ex}#2}="b",(7,2.2)*+{\;#1}="c"
\ar@{->} "a";"b" \ar @{->}"a";"c" \endxy }

\newcommand{\br}{\mbox{\rule{0.7em}{0.2ex}\hspace{-0.04em}\rule{0.08em}{1.7ex}}}

\newcommand{\tl}{\mbox{\rule{0.08em}{1.7ex}\rule[1.54ex]{0.7em}{0.2ex}}}

\newcommand{\hh}{\mbox{\rule{0.7em}{0.2ex}\hspace{-0.7em}\rule[1.5ex]{0.70em}{0.2ex}}}
\newcommand{\vv}{\mbox{\rule{0.08em}{1.7ex}\hspace{0.6em}\rule{0.08em}{1.7ex}}}

\newcommand{\sq}{\mbox{\rule{0.08em}{1.7ex}\hspace{-0.00em}\rule{0.7em}{0.2ex}\hspace{-0.7em}\rule[1.54ex]{0.7em}{0.2ex}\hspace{-0.03em}\rule{0.08em}{1.7ex}}}

\def\e{{e \rule[-0.5ex]{0em}{2ex}}}
\newcommand{\sqbox}{\hfill$\Box$}
\def\txtstyle{\def\labelstyle{\textstyle}}

\begin{document}
\title{Thin Elements and Commutative Shells \\ in Cubical \w-categories}
\author{Philip J. Higgins\thanks{email: p.j.higgins@durham.ac.uk}, \\
Department of Mathematical Sciences, \\ Science Laboratories, \\ South Rd., \\
Durham, DH1 3LE,  U.K.}

\maketitle

\begin{center}{\bf UWB Maths Preprint 04.16}
\end{center} 

\begin{abstract}The relationships between thin elements,
commutative shells and connections in cubical \w-categories are
explored by a method which does not involve the use of pasting
theory or nerves of \w-categories (both of which were previously
needed for this purpose; see \cite{A-B-S02}, Section 9). It is
shown that composites of commutative shells are commutative and
that thin structures are equivalent to appropriate sets of
connections; this work extends to all dimensions the results
proved in dimensions 2 and 3 in \cite{B-M99,B-M04}.
\end{abstract}

\section*{Introduction}

\noindent Thin structures in simplicial sets were introduced by
Dakin in \cite{D77} and were applied to cubical sets in
\cite{B-H77,B-H81,B-H81a}. In the cubical case a thin structure is
equivalent to an \w-groupoid structure.

In this paper we use the term \emph{thin structure} in a weaker
sense which is appropriate for cubical \w-categories and does not
imply the existence of inverses. This concept was introduced in
the 2-dimensional case in \cite{S-W}, as arising from a pair of
connections in dimension 2, and the equivalence of these notions
in this dimension was shown in \cite{B-M99}. We extend this result
to all dimensions (Theorem 3.1).

The definition of thin structure depends on the notion of
\emph{commutative cube} or \emph{commutative shell}. This was
studied in the 3-dimensional case in \cite{S-W} for certain
`special' double categories, and in \cite{B-M04} as part of the
proof of a van Kampen theorem for homotopy double groupoids. A key
result there (see also \cite[Proposition 3.11]{S-W})  was that any
composition of commutative 3-shells is commutative; we prove this
result in all dimensions. We also prove in the general case that
thin elements are precisely those which can be expressed as
composites using only identity elements and connections (Theorem
2.8).

These results, which were proved for \w-groupoids in \cite{B-H81},
can be deduced for \w-categories from results in \cite{A-B-S02}.
However, the methods used in \cite{A-B-S02} depend on the use of
pasting theory and nerves of \w-categories which tend to obscure
the intuitive nature of thinness and commutativity. The approach
used below is simpler and more direct. The basic simplification is
the use of a ``partial folding operation" $\Psi$ in place of the
full folding operation $\Phi_n$ used in \cite{A-B-S02}. The
operation $\Phi_n$ is needed to prove the equivalence between
cubical \w-categories with connections and globular \w-categories,
but the simpler $\Psi$ suffices for a detailed study of thinness
and commutativity. This direct approach should facilitate
applications to homotopy theory
(cf.\cite{B-M99,B-M04,Spencer,S-W}) and to concurrency theory in
computer science (cf. \cite{Gaucher,Goubault}).

\section{Composing the faces of a cube} \label{S:1}
Let $\C$ be a cubical \w-category as defined in \cite{A89} and
\cite{A-B-S02}; for the moment we do not assume the existence of
connections. If $x \in \C_n$ is an $n$-cube in $\C$ one may ask
which of its $(n-1)$-faces have common $(n-2)$-faces and can be
composed in $\C_{n-1}$. The answer is that the following pairs of
faces (and in general only these pairs) can be composed:
$$ (\partial^{-}_i x, \partial^{+}_{i+1} x), \qquad (\partial^{+}_i x, \partial^{-}_{i+1}
x), \qquad i=1,2,\ldots,n-1.$$ Thus the faces of $x$ (by which we
mean its $(n-1)$-faces) divide naturally into two sequences
$$ (\partial^{-}_1 x,\partial^{+}_2 x,\partial^{-}_3 x,\ldots,\partial^{\pm}_n
x) \mbox{ and } (\partial^{+}_1 x,\partial^{-}_2 x,\partial^{+}_3
x,\ldots,\partial^{\mp}_n x)$$ in which neighbouring pairs can be
composed. We call these respectively the \emph{negative} and the
\emph{positive} faces of $x$. This agrees with the terminology of
\cite{A-B-S02}.

We will frequently use 2-dimensional arrays of elements of $\C_n$,
and these will be shown in tabular form such as
$$ \vcenter
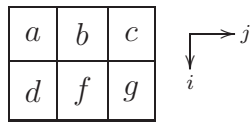
{\xymatrix@M=0pt @=0.5pc{ \ar @{-} [dddd] \ar @{-} [rrrrrr] &&\ar @{-}
[dddd] && \ar @{-} [dddd] && \ar @{-} [dddd] \\
& a & & b & & c & \\
\ar @{-} [rrrrrr] &&&&&& \\
& d& & f & & g &  \\
\ar @{-} [rrrrrr]  &&&&&&}}  \quad \directs{j}{i}$$ When using
this notation we will always assume that pairs of adjacent
elements in the table have a common face in the relevant
direction, in this case $i$ or $j$. Thus in the array above we
assume that $\partial ^+_ja = \partial ^-_j b$, $\partial ^+_ia =
\partial ^-_i d$, etc. The array can then be composed {\it by
rows}: $(a \circ _j b  \circ _j c)  \circ _i ( d  \circ _j f \circ
_jg)$, or {\it by columns}: $(a  \circ _i d)  \circ _j(b \circ _i
f) \circ _j(c  \circ _ig)$. It is a consequence of the interchange
law that these two elements of $\C_n$ are equal, and we call their
common value the {\it composite} of the array. To emphasise our
implicit assumptions, we will use the term {\it composable array}.

More general composites in the form of rectangular partitions of a
rectangle will also be used. The simplest is of the form
$$\vcenter
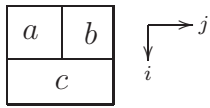
{\xymatrix@M=0pt @=0.5pc{\ar @{-} [dddd] \ar @{-} [rrrr] &&\ar
@{-} [dd] && \ar @{-} [dddd]\\
&a&&b& \\
\ar @{-} [rrrr]&&&&\\
&& c&& \\
\ar @{-} [rrrr]&&&&}} \quad   \directs{j}{i}$$ Here the implicit
assumptions are that $\partial ^+_ja= \partial ^-_j b$ and
$\partial ^+_i(a \circ _j b) = \partial ^-_i c$, and so the
composite $(a \circ _j b) \circ_i c$ can be formed. For a more
general (finite) partition of a rectangle by rectangles, labelled
by members of $\C_n$, we will assume that two elements, or
composites of elements, with a common edge in the partition, have
a  common face at this edge, so allowing the  composition in the
corresponding direction. For example, in the partition
$$\vcenter
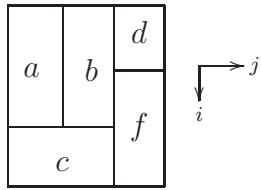
{\xymatrix@M=0pt @=0.5pc{\ar @{-} [dddddd] \ar @{-} [rrrrrr] &&\ar
@{-} [dddd] && \ar @{-} [dddddd]&& \ar @{-} [dddddd]\\
&&&& &d& \\
&a&&b&\ar @{-} [rr]&& \\
&&&&&&\\
\ar @{-} [rrrr]&&&&&f& \\
&& c&& &&\\
\ar @{-} [rrrrrr]&&&&&&}}  \quad { \directs{j}{i}}$$
 we assume, in addition to the relations above, that $$\partial^+_i
 d =\partial^-_i f \quad \text{ and } \quad \partial^+_j[(a \circ_j b)\circ_i c] =
 \partial ^-_j(d \circ_i f)$$ and so the composite $ [(a \circ_j b)\circ_i
 c]\circ_j(d \circ_i f)$ can be formed.  We shall call  such diagrams
 {\it composable partitions} provided that some sequence of
 compositions exist which combines all the elements to form a
 composite of the partition. We will specify the sequence if there
 is any ambiguity.

Now suppose that $\C$ is a \emph{cubical \w-category with
connections}, that is, it has for each $n$ and for
$i=1,2,\ldots,n$, extra structure maps $\Gamma^{+}_i,\Gamma^{-}_i:
\C_n \to \C_{n+1}$ (called connections) satisfying the identities
set out, for example,  in \cite{A89} and in Section 2 of
\cite{A-B-S02}. We shall make free use of the defining identities
in \cite{A-B-S02} without further comment. (In fact, the existence
of connections in $\C$
 implies that all composites of a given composable partition are
 equal. This is because, using connections, any composable
 partition can be refined to a composable {\it array}, and any
 composite of the partition must be equal to the unique composite
 of this array: see \cite{B-M99} and \cite{D-P}. We do not need
 to use this general theorem.)

The `degenerate' elements $\Gamma^{+}_{k} x, \Gamma^{-}_{k} x$ and
$\varepsilon_k x$ will sometimes be represented in composable
arrays by the symbols $\tl, \br$ and $(\vv \mbox{ or } \hh)$
respectively. The symbol $\sq \;$ will be used to denote an
element that is an identity for both the horizontal and vertical
compositions. These symbols will  be used only where the elements
of $\C_n$ they represent are uniquely determined by the
composability of the array. (The lines in these symbols are
designed to indicate that the corresponding faces are identities
in the direction of these edges.)

For example, the array
$$\vcenter{\xymatrix@M=0pt @=0.7pc{ \ar @{-} [dddd] \ar @{-} [rrrr]
&&\ar @{-} [dddd] && \ar @{-} [dddd] \\& x & & \hh & \\ \ar @{-}
[rrrr] &&&&
\\ & \vv& & \sq &  \\ \ar @{-} [rrrr]  &&&&}}
\raisebox{2ex}{\directs{j}{i}} \quad  \raisebox{0ex}{ denotes }
\quad \vcenter{\xymatrix@M=0pt @=0.5pc{ \ar @{-} [dddd] \ar @{-}
[rrrr] &&\ar @{-} [dddd] && \ar @{-} [dddd]
\\& x & & \varepsilon_j \partial^{+}_j x& \\ \ar @{-} [rrrr] &&&&
\\ & \varepsilon_i \partial^{+}_i x& & \varepsilon_i
\partial^{+}_j \varepsilon_j \partial^{+}_j x &  \\
\ar @{-} [rrrr]  &&&&}} \quad  {\directs{j}{i}} $$ whose composite
is $x$ itself; and the array
$$\def\labelstyle{\textstyle} \vcenter{\xymatrix@M=0pt @=0.7pc{ \ar @{-} [dddd] \ar @{-} [rrrr] &&\ar @{-} [dddd] && \ar @{-}
[dddd] |(0.25){a \rule[-0.5ex]{0em}{2ex}}  |(0.75){b
\rule[-0.3ex]{0em}{2ex}} \\& \tl & & \hh&   \\ \ar @{-} [rrrr]
&&&&&
\\ & \vv& & \tl &  \\ \ar @{-} [rrrr]|(0.25){ \rule{0.2em}{0ex}a}   |(0.75){\,b\,}   &&&&}}
{ \directs{i+1}{i}} \quad \text{ denotes \quad}
\vcenter{\xymatrix@M=0pt @=0.5pc{ \ar @{-} [dddd] \ar @{-} [rrrr]
&&\ar @{-} [dddd] && \ar @{-} [dddd]
\\& \Gamma^{+}_i a & & \varepsilon_{i+1} a& \\ \ar @{-} [rrrr] &&&&
\\ & \varepsilon_i a& & \Gamma^{+}_i b &  \\
\ar @{-} [rrrr]  &&&&}}\quad  \directs{i+1}{i} $$ whose composite,
according to the transport law \cite[(2.6)]{A-B-S02}, is
$\Gamma^{+}_i(a \, \o_i \, b)$. (Here the labels on the edges
denote the corresponding faces of elements of the array.)

The following  identities (see \cite[(2.7)]{A-B-S02}) will be
important in what follows:
$$
\txtstyle \vcenter{\xymatrix@M=0pt @=0.7pc{ \ar @{-} [dd] \ar @{-}
[rrrr]|(0.75){\,a\,} && \ar @{-} [dd] && \ar @{-} [dd] \\& \tl & & \br&   \\
\ar @{-} [rrrr] |(0.25){\,a\,} &&&&}} \; = \;
\vcenter{\xymatrix@M=0pt @=0.7pc{ \ar @{-} [dd]
\ar @{-} [rr]|{\;a\,}  && \ar @{-} [dd] \\& \vv &\\
\ar @{-} [rr] | {\;a\,}&&}} \raisebox{1mm}{ \directs{i+1}{i}}
\raisebox{0mm}{ and \quad} \vcenter{\xymatrix@M=0pt @=0.7pc{ \ar
@{-} [dddd]|(0.75){ a \rule[-0.5ex]{0em}{2ex}}  \ar @{-} [rr] &&
\ar @{-}
[dddd] |(0.25){a \rule[-0.5ex]{0em}{2ex}}\\& \tl &    \\ \ar @{-} [rr] && \\
a & \br & \\
\ar @{-} [rr] &&}} \; = \; \vcenter{\xymatrix@M=0pt @=0.7pc{ \ar
@{-} [dd]|{a \rule[-0.5ex]{0em}{2ex}}
\ar @{-} [rr] && \ar @{-} [dd]|{a \rule[-0.5ex]{0em}{2ex}}  \\& \hh &\\
\ar @{-} [rr] &&}} \raisebox{1mm}{ \directs{i+1}{i}}
$$
that is, $\Gamma^{+}_i a \; \o_{i+1} \; \Gamma^{-}_i a =
\varepsilon_i a$ and $\Gamma^{+}_i a \; \o_i \; \Gamma^{-}_i a =
\varepsilon_{i+1}a$.

 \vspace{2ex}

 We now define the \emph{elementary
folding operations} $\psi_i: \C_n \to \C_n, \,(i=1,2,\ldots,n-1)$
by
$$ \psi_i x =
\vcenter{\xymatrix@M=0pt @=0.7pc{ \ar @{-} [dd] \ar @{-} [rrrrrr]
&&\ar
@{-}[dd] && \ar @{-} [dd] && \ar @{-} [dd] \\
& \tl & & x & & \br & \\
\ar @{-} [rrrrrr] &&&&&&}} \quad   \directs{i+1}{i}
$$
The chief identities satisfied by these operations are set out in
Proposition 3.3(i) of \cite{A-B-S02}; we shall not need the
``braid relations" proved in Theorem 5.2 of that paper. As the
picture suggests, the effect of $\psi_i$ is to ``fold" the faces
$\partial^{-}_{i+1} x$ and $\partial^{+}_{i+1} x$ to the $i^{th}$
direction so that they abut the faces $\partial^{+}_i x$ and
$\partial^{-}_i x$ respectively, and to compose the two pairs of
faces. The composition of all the negative (and positive) faces of
$x$, together with certain faces of the connections,  can
therefore be achieved by the folding operation $\Psi: \C_n
\to\C_n$ defined by $\Psi = \psi_1 \psi_2 \ldots \psi_{n-1}$.

We note that $\Psi$ is not the full folding operation $\Phi_n$ of
\cite{A-B-S02} (which maps $\C_n$ into its globular part (see
Proposition 3.5 of \cite{A-B-S02}), but $\Psi$ is sufficient for
the study of thin elements and commutative shells. We will return
to this point later. For now, the two most important faces of
$\Psi x$, namely
$$Px  = \partial^{+}_1 \Psi x \mbox{\; and \;} Nx  = \partial^{-}_1 \Psi x,$$
are to be viewed as convenient embodiments of the positive and
negative boundaries of $x$.

\begin{lemma} \label{L:1.1}
\mbox{ }
\begin{enumerate} [\rm (i)]
\item $\psi_1 \psi_2 \ldots \psi_{r-1} \varepsilon_r =
\varepsilon_1: \C_{n-1} \to \C_n \mbox{ for } 1 \leqslant r
\leqslant n-1.$ \item If $y \in \C_{n-1}$ then $\Psi \varepsilon_j
y \in \varepsilon_1 \C_{n-1} \mbox{ for } 1 \leqslant j \leqslant
n-1$.
\end{enumerate}
\end{lemma}

\Proof (i) This follows immediately from the identities $\psi_i
\varepsilon_{i+1} = \varepsilon_i$.\\
(ii) similarly, using $\psi_i \varepsilon_j = \varepsilon_j
\psi_{i-1}$ for $j < i$ and $\psi_j \varepsilon_j =
\varepsilon_j$, we have
\begin{align*}
\Psi \varepsilon_j&=  \psi_1 \psi_2 \ldots \psi_{n-1}
\varepsilon_j \\
&=   \psi_1 \psi_2 \ldots \psi_j \varepsilon_j \psi_j \psi_{j+1}
\ldots \psi_{n-2}\\
&=   \psi_1 \psi_2 \ldots \psi_{j-1} \varepsilon_j \psi_j
\psi_{j+1} \ldots \psi_{n-2}\\
&=   \varepsilon_j \psi_j \psi_{j+1} \ldots \psi_{n-2} \mbox{ by
(i).} \tag*{$\Box$}
\end{align*}

If $x \in \C_n$, the \emph{shell} $\bpartial x$ of $x$ is the
family  consisting of all its faces $\partial^{\alpha}_i x \;
(i=1,2,\ldots,n; \alpha = + , -)$.

\begin{prop}
Let $\C$ be a cubical \w-category with connections and $x \in
\C_n$. Then
\begin{enumerate} [\rm (i)]
\item all faces $\partial^{\alpha}_i \Psi x$ with $i \geqslant 2$
are of the form $\varepsilon_1 z^{\alpha}_i$, where $z^{\alpha}_i
\in \C_{n-2}$; \item $\bpartial Nx  = \bpartial Px $; \item
$\bpartial \Psi x$ is uniquely determined by $Nx $ and $Px $.
\end{enumerate}
\end{prop}

\Proof
\begin{enumerate} [\hspace{-1em} (i)]
\item We have $\partial^{\alpha}_i \psi_j = \psi_j
\partial^{\alpha}_i$ for $i > j+1$, so for $i \geqslant 2$,
$$ \partial^{\alpha}_i \Psi x = \partial^{\alpha}_i \psi_1 \ldots
\psi_{n-1} x = \psi_1 \psi_2 \ldots \psi_{i-2} \partial^{\alpha}_i
\psi_{i-1}, \ldots, \psi_{n-1} x,$$ \noindent (if $i =2$, the
$\psi_1 \ldots, \psi_{i-2}$ are missing). But $\partial^{\alpha}_i
\psi_{i-1} = \varepsilon_{i-1}
\partial^{\alpha}_{i-1} \partial^{\alpha}_i$, so
\begin{align*}
\partial^{\alpha}_i \Psi x &= \psi_1 \psi_2 \ldots \psi_{i-2}
\varepsilon_{i-1} z^{\alpha}_i \mbox{ where } z^{\alpha}_i \in
\C_{n-2} \\
&= \varepsilon_1 z^{\alpha}_i \mbox{ by Lemma1.1 (i) } .
\end{align*}
\item This follows from (i) because
$$ \partial^{\alpha}_i Nx  = \partial^{\alpha}_i \partial^{-}_1
\Psi x = \partial^{-}_1 \partial^{\alpha}_{i+1} \Psi x =
\partial^{-}_1 \varepsilon_1 z^{\alpha}_{i+1} = z^{\alpha}_{i+1}$$
and similarly $\partial^{\alpha}_i P_x = z^{\alpha}_{i+1}$. \item
The faces of $\Psi x$ are $Nx,\,  Px $ and the elements
$\varepsilon_1 z$ where $z$ is a face of $Nx$  (or $Px$). \sqbox\\
\end{enumerate}

 In the abstract, an $\emph{$n$-shell}$ in $\C$ is a
family $\s = \{ s^{\alpha}_i; s^{\alpha}_i \in \C_{n-1}; i =
1,2,\ldots,n; \alpha = +, -\}$ where the $s^{\alpha}_i$ satisfy
the incidence relation
$$ \partial^{\beta}_j s^{\alpha}_i = \partial^{\alpha}_{i-1}
s^{\beta}_j \mbox{ for } 1 \leqslant j < i \leqslant n \mbox{ and
} \alpha, \beta \in \{+, -\}.$$ We denote by $\square \,\C_{n-1}$
the set of such shells. The usual cubical incidence relations
imply that $\bpartial x \in \square \,\C_{n-1}$ for all $x \in
\C_n$.

If $\{\C_0, \C_1, \ldots, \C_{n-1}\}$ is a cubical
$(n-1)$-category with connections, then $\square \C_{n-1}$ has
naturally defined operations $\o_i \;(1 \leqslant i \leqslant n)$,
and connections $\mathbf{\Gamma}^{\alpha}_j:\C_{n-1} \to \square
\C_{n-1} \; (1 \leqslant j \leqslant n-1)$ which, together with
the obvious structure maps $\partial^{\alpha}_i: \square \,
\C_{n-1} \to \C_{n-1}$ and $\bvarepsilon_j: \C_{n-1} \to \square
\, \C_{n-1}$, make $\{\C_0, \C_1,\ldots, \C_{n-1}, \square \,
\C_{n-1}\}$ a cubical $n$-category with connections (cf.
\cite{B-H81}, section 5). Thus we can define folding maps $\psi_i,
\Psi: \square \C_{n-1} \to \square \C_{n-1}$ which obviously
satisfy:

\begin{lemma} \label{L:1.2}
In a cubical $n$-category $(\C_1, \C_2, \ldots, \C_n)$ with
connections, the map $$x \mapsto \bpartial x : \C_n \to \square
\C_{n-1},$$ together with identity maps in lower dimensions, gives
a morphism of cubical $n$-categories with connections from $(\C_1,
\C_2, \ldots, \C_n)$ to  $(\C_1, \C_2, \ldots, \C_{n-1}, \square
\C_{n-1})$. In particular
$$ \mathbf{\Gamma}^{\alpha}_i x = \bpartial \Gamma^{\alpha}_i x,
\;\psi_j \bpartial x = \bpartial \psi_j x, \;\Psi \bpartial x =
\bpartial \Psi x, $$
$$ N \bpartial x = \bpartial N x, \; P \bpartial x = \bpartial Px \mbox{ and }
\bpartial \varepsilon_j x = \bvarepsilon_j x.$$ \sqbox
\end{lemma}
\begin{theorem} \label{T:1.4}
Let $\C$ be a cubical \w-category (or  a cubical $m$-category)
with connections. Let  $a \in \C_n$ and $\s \in \square \C_{n-1}$.
A necessary and sufficient condition for the existence of $x \in
\C_n$ such that
$$ \bpartial x = \s \mbox{ and } \Psi x = a $$
is that
$$ \Psi \s = \bpartial a.$$
If $x$ exists, it is unique.
\end{theorem}

\Proof The necessity of $\Psi \s = \bpartial a$ follows from
$\bpartial \Psi x = \Psi \bpartial x$ (Lemma 1.2). The existence
and uniqueness will be deduced from:

\begin{lemma} \label{L:1.5}
Let $a \in \C_n$ and $\s \in \square \C_{n-1}$ satisfy $\bpartial
a = \psi_j \s$ for some $j \in \{1,2,\ldots,n-1\}$. Then there is
a unique $x$ in $\C_n$ such that $\bpartial x = \s$ and $\psi_j x
= a$.
\end{lemma}

\Proof First suppose that $x$ exists, and consider the  array
$$
A: \quad \vcenter{\xymatrix@M=0pt @=0.7pc{ \ar @{-} [dddddd] \ar
@{-} [rrrrrr] &&\ar @{-}
[dddddd] && \ar @{-} [dddddd] && \ar @{-} [dddddd] \\
& \sq & & \vv & & \tl & \\
\ar @{-} [rrrrrr] &&&&&& \\
& \tl& & x & & \br &  \\
\ar @{-} [rrrrrr]  &&&&&& \\
& \br& & \vv & & \sq &  \\
\ar @{-} [rrrrrr]  &&&&&&}} \quad  \directs{j+1}{j}
$$
where the elements surrounding $x$ are determined by the faces of
$x$ so that all rows and columns are composable. The composite of
the middle row is $\psi_j x = a$. The elements of the first and
third rows are determined by the faces of $x$, i.e. by $\s$. Hence
the composite of $A$ is determined by $a$ and $\s$. But if we
compose $A$ by columns and use the law $\mathbf{\Gamma}^{+}_j t \,
\o_j \, \mathbf{\Gamma}^{-}_j t = \bvarepsilon_{j+1} t$, we see
that the composite of $A$ is $x$ itself. Hence $x$ is unique.

To prove existence, we note that the array $A$ gives a formula for
$x$ in terms of $a$ and $\s$.

\noindent So, given $a$ and $\s$ we define $x$ to be the composite
of the composable partition
$$
x= \quad \vcenter{\xymatrix@M=0pt @=0.5pc{ \ar @{-}
[dddddd] \ar @{-} [rrrr] &&\ar @{-}[dd] && \ar @{-} [dddddd] \\
& \varepsilon_j s^{-}_j & & \Gamma^{+}_j s^{+}_{j+1} & \\
\ar @{-} [rrrr] &&&& \\
&& a \rule[-0.5ex]{0mm}{2ex} && \\
\ar @{-} [rrrr]  && \ar @{-}[dd] && \\
& \Gamma^{-}_j s^{-}_{j+1}& & \varepsilon_j s^{+}_j & \\
\ar @{-} [rrrr]  &&&&}} \quad { \directs{j+1}{j}}
$$
Here the first and third rows are the same as those in the array
$A$, except that the 2-fold identities $\sq$ have been omitted,
being redundant. Because we are assuming that $\bpartial a= \psi_j
\s$, the faces $\partial ^-_j a$ and $\partial ^+_j a$ are the
same as the upper and lower faces of the composite $\vcenter{
\xymatrix@M=0pt @=0.3pc{\ar @{-} [rrrrrr] \ar@{-} [dd]  && \ar@{-} [dd]&&\ar@{-} [dd]& &\ar@{-} [dd] \\
& \tl && x && \br &\\
\ar @{-} [rrrrrr] &&&&&&}}$ in $A$. Hence the partition is
composable and we can compute $x$ from $a$ and $\s$, by composing
rows first. It remains to verify that $\psi_j x = a$ and
$\bpartial x = \s$.

The faces $\partial^{\pm}_{j+1} x$ of $x$ are $\o_j$-composites as
indicated in the diagram
$$
x= \quad\txtstyle \vcenter{\xymatrix@M=0pt @=0.7pc{ \ar @{-}
[dddddd] |(0.2)\e |(0.5)\e |(0.8){s^{+}_{j+1}
\rule[-1ex]{0mm}{3ex} } \ar @{-} [rrrr] &&\ar @{-}[dd] &&\ar @{-}
[dddddd]
 |(0.2){s^{+}_{j+1} \rule[-1ex]{0mm}{2.5ex}  }|(0.5)\e |(0.8)\e\\
 & \rule{0mm}{3ex} \qquad & & \qquad & \\
\ar @{-} [rrrr] &&&& \\
 &&  a && \\
\ar @{-} [rrrr]  && \ar @{-}[dd] && \\
& \rule{0mm}{3ex} \qquad & & \qquad & \\
\ar @{-} [rrrr]  &&&&}} \quad { \directs{j+1}{j}}
$$
where the faces labelled $e$ are identities for $\o_j$. (Note that
$\partial^{\pm}_{j+1} a$ are identities because $\bpartial a =
\psi_j \s)$. Hence
$$
\psi_j x = \quad \vcenter{ \xymatrix@M=0pt @=0.5pc{\ar @{-}
[dddddd] \ar @{-} [rrrrrrrr] &&\ar @{-}[dddddd] && \ar@{-} [dd] &&
\ar @{-} [dddddd] &&
\ar @{-} [dddddd]\\
& \sq && \varepsilon_j s^{-}_j && \Gamma^{+}_j s^{+}_{j+1} && \Gamma^{-}_j s^{+}_{j+1} & \\
\ar @{-} [rrrrrrrr] &&&&&&&& \\
& \sq &&& a \rule[-1ex]{0em}{3ex}  &&& \sq & \\
\ar @{-} [rrrrrrrr] &&&&\ar@{-}
[dd] &&&& \\
& \Gamma^{+}_j s^{-}_{j+1} && \Gamma^{-}_j s^{-}_{j+1} && \varepsilon_j s^{+}_j && \sq & \\
\ar @{-} [rrrrrrrr]  &&&&&&&&}} \quad { \directs{j+1}{j}}
$$
Here, the first and third columns are the relevant connections
expanded by the transport law and then simplified. The diagram can
be viewed as a $3 \times 3$ array in which two elements happen to
be horizontal composites; therefore we may compute $\psi_j x$ by
composing the rows first instead of the columns.  Using the law
$\Gamma^{+}_j t \, \o_{j+1} \, \Gamma^{-}_j t = \varepsilon_j t$,
we find that $\psi_j x = a$, as required.

Finally $\bpartial x$ is the composite (by rows) of
$$
\vcenter{\xymatrix@M=0pt @=0.5pc{ \ar @{-}
[dddddd] \ar @{-} [rrrr] &&\ar @{-}[dd] && \ar @{-} [dddddd] \\
& \bpartial \varepsilon_j s^{-}_j & & \bpartial \Gamma^{+}_j s^{+}_{j+1} & \\
\ar @{-} [rrrr] &&&& \\
&& \rule[-1ex]{0mm}{3ex} \bpartial a \rule[-1ex]{0mm}{3ex}&& \\
\ar @{-} [rrrr]  && \ar @{-}[dd] && \\
& \bpartial \Gamma^{-}_j s^{-}_{j+1} & & \bpartial \varepsilon_j
s^{+}_j & \\
\ar @{-} [rrrr]  &&&&}} \quad  \directs{j+1}{j}
$$
Since $\bpartial a = \psi_j \s$, this is the composite in $\square
\, \C_{n-1}$ of the array
$$
\vcenter{\xymatrix@M=0pt @=0.7pc{ \ar @{-} [dddddd] \ar @{-}
[rrrrrr] &&\ar @{-}
[dddddd] && \ar @{-} [dddddd] && \ar @{-} [dddddd] \\
& \sq & & \vv & & \tl & \\
\ar @{-} [rrrrrr] &&&&&& \\
& \tl& & \mathsf{s} & & \br &  \\
\ar @{-} [rrrrrr]  &&&&&& \\
& \br& & \vv & & \sq &  \\
\ar @{-} [rrrrrr]  &&&&&&}} \quad  \directs{j+1}{j}
$$
which, as before, is just $\s$. \sqbox\\

We now use induction on $r \leqslant n-1$ to show that if
$\bpartial a = \psi_1 \psi_2 \ldots \psi_r \s$, then there is a
unique $x \in \C_n$ such that $\psi_1 \psi_2 \ldots \psi_r x = a$
and $\bpartial x = \s$. The case $r=1$ is covered by  Lemma 1.5.
Suppose that the result is true when $r = t-1 < n-1$ and that
$\bpartial a = \psi_1 \psi_2 \ldots \psi_r \s$. Then, by induction
hypothesis, there is a unique $y \in \C_n$ such that $\psi_1
\psi_2 \ldots \psi_{r-1} y = a$ and $\bpartial y = \psi_t \s$.
But, again by Lemma 1.5, there is then a unique $x \in \C_n$ with
$\psi_t x = y$ and $\bpartial x = \s$, completing the induction.
The case $r=n-1$, completes the proof of the theorem. \sqbox

\section{Thin elements and commutative shells}

We say that an element $x \in \C_n$ is {\it thin} if $\Psi x \in
\varepsilon_1 \C_{n-1}$. Then $\Psi x = \varepsilon_1 Nx  =
\varepsilon_1 Px $.

We say that a shell $\s \in \square \, \C_{n-1}$ is
\emph{commutative} if $N \s = P \s$.

\begin{prop} \label{P:2.1}
Let $\C$ be a cubical $\omega$-category  (or cubical $m$-category)
with connections.
\begin{enumerate} [\rm (i)]
\item The shell of a thin element of $\C_n$ is a commutative
$n$-shell. \item A commutative $n$-shell is the same thing as a
thin element of $\square \, \C_{n-1}$. \item Any commutative
$n$-shell $\s$ has a unique thin filler (i.e. a thin element $x
\in \C_n$ with $\bpartial x = \s$).
\end{enumerate}
\end{prop}

\Proof (i) Let $x \in \C_n$ be thin. Then $\Psi x = \varepsilon_1 z$ for some
$z \in \C_{n-1}$, so $Nx=Px =z$. Therefore $N \bpartial x= \bpartial Nx
= \bpartial Px = P \bpartial x$, by Lemma 1.3.   \\
(ii) If $\s \in \square \, \C_{n-1}$ is commutative, then $N \s =
P \s = u$, say. The other faces of $\Psi \s$ are all of the form
$\varepsilon_1 v$, where $v$ is a face of $u$. These faces
determine the shell $\Psi \s$ and identify it as $\bvarepsilon_1
u$, so $\s$ is thin. The converse is obvious.\\
(iii) Let $\s$ be a commutative $n$-shell. Let $u = N \s = P \s$
and put $a = \varepsilon_1 u \in \C_n$. Then $\bpartial a =
\bvarepsilon_1 u = \Psi \s$, by (ii) and Lemma 1.3. By Theorem
1.4, there is a unique $x \in \C_n$ with $\bpartial x = \s$ and
$\Psi x = a = \varepsilon_1 u$. \sqbox

\begin{prop} \label{P:2.2}
Let $\C$ be  a cubical \w-category (or a cubical $m$-category)
with connections.
\begin{enumerate} [\rm (i)]
\item Elements in $\C$ of the form $\varepsilon_i c$ or
$\Gamma^{\alpha}_i c$ are thin. \item If $a, b \in \C_n$ are thin
and $c = a \circ_{i} b$ then $c$ is thin.
\end{enumerate}
\end{prop}

\Proof (i) The thinness of $\varepsilon_i c$ has been proved in
Lemma 1.1. To prove thinness of $\Gamma^{\alpha}_i c$, we first
establish some formulae involving the $\psi_i$ and
$\Gamma^{\alpha}_i$.

\begin{lemma} \label{L:2.3}
\begin{enumerate} [\rm (i)]
\item $\psi_i \Gamma^{\alpha}_i = \varepsilon_i$. \item $\psi_j
\Gamma^{\alpha}_i = \Gamma^{\alpha}_i \psi_{j-1}$, if $j > i+1$.
\item $\psi_i \psi_{i+1} \Gamma^{+}_i c =
\varepsilon_i(\Gamma^{+}_i \partial^{-}_{i+1} c \; \o_{i+1} \;
c)$.\\
$\psi_i \psi_{i+1} \Gamma^{-}_i c = \varepsilon_i(c \; \o_{i+1} \;
\Gamma^{-}_i \partial^{+}_{i+1} c).$
\end{enumerate}
\end{lemma}

\Proof (i) \vspace{-3.5ex}
\begin{align*} \psi_i \Gamma^{+}_i c&=  \Gamma^{+}_i
\partial^{-}_{i+1} \Gamma^{+}_i c \; \o_{i+1} \; \Gamma^{+}_i c \;
\o_{i+1} \;
\Gamma^{-}_i \partial^{+}_{i+1} \Gamma^{+}_i c. \\
&=   \Gamma^{+}_i \varepsilon_i \partial^{-}_{\alpha} c \;
\o_{i+1}
\; \Gamma^{+}_i c \; \o_{i+1} \; \Gamma^{-}_i c \\
&=   \varepsilon_{i+1} \varepsilon_{i} \partial^{-}_{\alpha} c \;
\o_{i+1} \; \Gamma^{+}_i c \; \o_{i+1} \; \Gamma^{-}_i c \\
&=   \Gamma^{+}_i c \; \o_{i+1} \; \Gamma^{-}_i c \\ &=
\varepsilon_i c.
\end{align*}
The proof for $\Gamma^{-}_i$ is similar.\\
(ii) This is proved similarly by using standard laws from pp.
80-81
of \cite{A-B-S02}.\\
(iii) When we try to compute $\psi_{i+1} \Gamma^{+}_i$, we are
hindered by the lack of a simple law involving $\Gamma^{-}_{i+1}
\Gamma^{+}_i c$. However, when we compute $\psi_i \psi_{i+1}
\Gamma^{+}_i c$, this difficulty disappears. We note that if $i+1
< j$ then $\psi_i ( a \; \o_j \; b) = \psi_i a \; \o_j \; \psi_i
b$. So
\begin{align*}
\psi_i \psi_{i+1} \Gamma^{+}_i c&=  \psi_i (\Gamma^{+}_{i+1}
\partial^{-}_{i+2} \Gamma^{+}_i c \; \o_{i+2} \; \Gamma^+_i c \;
\o_{i+2} \; \Gamma^{-}_{i+1} \partial^{+}_{i+2}
\Gamma^{+}_i c)\\
&= \psi_i \Gamma^{+}_{i+1} \partial^{-}_{i+2}\Gamma^{+}_i c \;
\o_{i+2} \; \psi_i \Gamma^+_i c \; \o_{i+2} \; \psi_i
\Gamma^{-}_{i+1}
\partial^{+}_{i+2} \Gamma^{+}_i c \\
&= \psi_i \Gamma^{+}_{i+1} \Gamma^{+}_i \partial^{-}_{i+1} c \;
\o_{i+2} \; \varepsilon_i c \; \o_{i+2} \; \psi_i \Gamma^{-}_{i+1}
\Gamma^{+}_i \partial^{+}_{i+1} c.
\end{align*}
We calculate the first and last terms separately:
\begin{align*}
\psi_i \Gamma^{+}_{i+1} \Gamma^{+}_i&=  \psi_i \Gamma^{+}_i
\Gamma^{+}_i = \varepsilon_i \Gamma^{+}_i \mbox{ by (i) } \\
\intertext{and }  \psi_i \Gamma^{-}_{i+1} \Gamma^{+}_i y&=
\Gamma^{+}_i
\partial^{-}_{i+1} \Gamma^{-}_{i+1} \Gamma^{+}_i y \; \o_{i+1} \;
\Gamma^{-}_{i+1} \Gamma^{+}_i y \; \o_{i+1} \; \Gamma^{-}_i
 \partial^{+}_{i+1}\Gamma^-_{i+1} \Gamma^{+}_i y \\
&=   (\Gamma^{+}_i \Gamma^{+}_i y \; \o_{i+1} \; \Gamma^{-}_{i+1}
\Gamma^{+}_i y ) \; \o_{i+1} \; \Gamma^{-}_i \varepsilon_{i+1}\partial^+_{i+1} \Gamma^+_i y \\
&=   (\Gamma^{+}_{i+1} \Gamma^{+}_i y \; \o_{i+1} \;
\Gamma^{-}_{i+1} \Gamma^{+}_i y ) \; \o_{i+1} \;  \Gamma^{-}_i\varepsilon_{i+1} y \\
&=   \varepsilon_{i+2} \Gamma^{+}_i y  \; \o_{i+1} \;
\varepsilon_{i+2} \Gamma^{-}_i y  \\
&=   \varepsilon_{i+2} (\Gamma^{+}_i y \; \o_{i+1} \; \Gamma^{-}_i
y) \\
&=   \varepsilon_{i+2} \varepsilon_i y.
 \intertext{Hence}
 \psi_i \psi_{i+1} \Gamma^{+}_i c&= ( \varepsilon_i \Gamma^{+}_i
\partial^{-}_{i+1} c \; \o_{i+2} \; \varepsilon_i c) \; \o_{i+2} \;
\varepsilon_{i+2} \varepsilon_i \partial ^+_{i+1} c \\
&=   \varepsilon_i(\Gamma^{+}_i \partial^{-}_{i+1} c \; \o_{i+1}
\; c).
\end{align*}
The proof is similar for $\Gamma^{-}_i$. This completes the proof
of Lemma 2.3. \hfill $\Box$ \\

Returning to Proposition 2.2(i), the proof that $\Gamma^{\alpha}_i
c$  (for $ c \in \C_{n-1})$ is thin is now straightforward:
\begin{align*}
\Psi \Gamma^{\alpha}_i c&=  \psi_1 \psi_2 \ldots \psi_{n-1}
\Gamma^{\alpha}_i c \\
&=   \psi_1 \psi_2 \ldots \psi_i \psi_{i+1} \Gamma^{\alpha}_i c
\psi_{i+1} \ldots \psi_{n-2} c, \mbox{ by 2.3(ii) }\\
&=   \psi_1 \psi_2 \ldots \psi_{i-1} \varepsilon_i z \mbox{ for
some $z \in \C_{n-1}$, by 2.3(iii)} \\
&=   \varepsilon_1 z \mbox{\quad by 1.1(i), as required. }
\end{align*}

\noindent (ii) To prove that composites of thin elements are thin
we introduce a subsidiary definition: if $ x \in \C_n$ and $1
\leqslant j \leqslant n-1$, we say that $x$ is \emph{$j$-thin} if
$\psi_1 \psi_2 \ldots \psi_j x \in \varepsilon_1 \C_{n-1}$. Thus,
for elements of $\C_n, \;(n-1)$-{\it thin} means \emph{thin}. We
take ``$x \in \C_n$ is $0$-thin" to mean $x \in \varepsilon_1
\C_{n-1}$.

\begin{lemma} \label{L:2.4}
If $j \geqslant 1$, then $x$ is $j$-thin if and only if $\psi_j x$
is $(j-1)$-thin. \sqbox
\end{lemma}
\begin{lemma} \label{L:2.5}
If $x=y  \o_i  z$ in $ \C_n$ and $y,z \in \varepsilon_j
\,\C_{n-1}$ then $x \in \varepsilon_j \, \C_{n-1}$.
\end{lemma}
(Note that in the case $i=j$, the hypotheses imply $x=y=z$.)\sqbox
\begin{lemma} \label{L:2.6}
If $x \in \varepsilon_k \, \C_{n-1}$, then $x$ is $(k-1)$-thin.
\end{lemma}

\Proof $\psi_1 \psi_2 \ldots \psi_{k-1} x = \psi_1 \psi_2 \ldots
\psi_{k-1} \varepsilon_k y = \varepsilon_1 y$ by (1.1)(i).
\sqbox\\

\noindent We now use induction on $j$ to prove:\\
\begin{equation*}
 \text{If }  a, b \in \C_n \text{ are }  j\text{-thin and } c= a
\o_i  b  \text{ for some }   1 \leqslant i \leqslant n \text{ then
} c \mbox{ is } j\text{-thin. } \tag{*}\end{equation*}

The case $j=0$ is contained in Lemma 2.5.

Suppose that (*) is true for $j=0,1,\ldots, k-1$, where $1
\leqslant k \leqslant n-1$. We will deduce (*) for $j=k$.

Let $c = a  \o_i  b$ in $\C_n$ and assume that $a$ and $b$ are
$k$-thin. We examine $\psi_k c = \psi_k(a \, \o_i \, b)$ in order
to prove that it is $(k-1)$-thin.\\

\noindent \textbf{Case 1:} $k < i-1$ or $k > i$. In this case
$\psi_k c = \psi_k a \, \o_i \, \psi_k b$ and $\psi_k a$ and
$\psi_k b$ are $(k-1)$-thin. By induction hypothesis, $\psi_k c$
is $(k-1)$-thin and so $c$ is $k$-thin by 2.4.\\

\noindent \textbf{Case 2:} $k=i-1$. Then
\begin{align*}
\psi_k c &= \;  \psi_{i-1} (a \, \o_i \, b) \\
 & \\
 & =\;  \vcenter{ \xymatrix@M=0pt @=0.7pc{ \ar @{-}
[dd] \ar @{-} [rrrrrrrr] && \ar @{-} [dd] && \ar @{-} [dd] &&
\ar @{-} [dd] && \ar @{-} [dd] \\
& \tl& & a & & b & & \br &  \\
\ar @{-} [rrrrrrrr]  &&&&&&&&&&}}  \quad  \directs{i}{i-1} \\
&\\ & =\; \vcenter{ \xymatrix@M=0pt @=0.7pc{ \ar @{-} [dddd] \ar
@{-} [rrrrrrrrrr] && \ar @{-} [dddd] && \ar @{-} [dddd] && \ar
@{-} [dddd]
&& \ar @{-} [dddd] && \ar @{-} [dddd] \\
& \sq  & & \vv  & & \tl & & b & & \br &  \\
\ar @{-} [rrrrrrrrrr]  &&&&&&&&&& \\
& \tl & & a & & \br & & \vv  & & \sq  &  \\
\ar @{-} [rrrrrrrrrr]  &&&&&&&&&& }} \quad  \directs{i}{i-1}\\
&\\
&=\;  \vcenter{ \xymatrix@M=0pt @=0.7pc{ \ar @{-} [dddd] \ar
@{-} [rrrrr] && \ar @{-} [dd] && & \ar @{-} [dddd] \\
& \varepsilon_{i-1} u &&&**[l] \psi_{i-1} b &  \\
\ar @{-} [rrrrr]  &&&\ar@{-}[dd]&& \\
&**[r] \psi_{i-1} a && &\varepsilon_{i-1} v &  \\
\ar @{-} [rrrrr]  &&&&& }} \quad  \directs{i}{i-1}\\
&\\
&= \; \vcenter{\xymatrix@M=0pt @=0.7pc{ \ar @{-} [dddd] \ar
@{-} [rrrrr] && \ar @{-} [dd] && & \ar @{-} [dddd] \\
& \varepsilon_{k} u &&&**[l] \psi_k b &  \\
\ar @{-} [rrrrr]  &&&\ar@{-}[dd]&& \\
&**[r] \psi_{k} a && &\varepsilon_{k} v &  \\
\ar @{-} [rrrrr]  &&&&& }} \quad  \directs{i}{i-1}
\end{align*}
where $u = \partial^{-}_{i-1} a, v = \partial^{+}_{i-1} b$. Now
$\psi_k a$ and $\psi_k b$ are $(k-1)$-thin, by 2.4, and
$\varepsilon_k u, \varepsilon_k v$ are $(k-1)$-thin, by 2.6. So
$\psi_k c$, being a composite of these, is $(k-1)$-thin by
induction hypothesis. Hence $c$ is $k$-thin.\\
\textbf{Case 3:} $k=i$. This is similar using the formula

\begin{align*}
\psi_k c &= \;  \psi_i (a \, \o_i \, b) \\
\quad\\
& =\; \vcenter{ \xymatrix@M=0pt @=0.5pc{ \ar @{-} [dddd] \ar @{-}
[rrrrrr] && \ar @{-} [dddd] && \ar @{-} [dddd] && \ar @{-}
[dddd]\\
& & & a & & &  \\
& \quad \tl \quad & \ar @{-} [rr]  &&& \quad \br \quad & \\
& & & b & & &  \\
\ar @{-} [rrrrrr]  &&&&&& }} \quad  \directs{i+1}{i}\\
&\\ & =\; \vcenter{ \xymatrix@M=0pt @=0.7pc{ \ar @{-} [dddd] \ar
@{-} [rrrrrrrrrr] && \ar @{-} [dddd] && \ar @{-} [dddd] && \ar
@{-} [dddd]
&& \ar @{-} [dddd] && \ar @{-} [dddd] \\
& \tl & & \hh & & a & & \br & & \vv  &  \\
\ar @{-} [rrrrrrrrrr]  &&&&&&&&&& \\
& \vv  & & \tl & & b & & \hh & & \br &  \\
\ar @{-} [rrrrrrrrrr]  &&&&&&&&&& }} \quad  \directs{i+1}{i}\\
&\\ & =\; \vcenter{ \xymatrix@M=0pt @=0.7pc{ \ar @{-} [dddd] \ar
@{-} [rrrrr] &&& \ar @{-} [dd] && \ar @{-} [dddd] \\
& **[r] \psi_i a &&& \varepsilon_i s &  \\
\ar @{-} [rrrrr] & & \ar @{-} [dd] &&&& \\
& \varepsilon_i t && &**[l] \psi_i b &  \\
\ar @{-} [rrrrr]  &&&&& }} \quad \quad  \directs{i+1}{i}
\end{align*}
where $s = \partial^{+}_{i+1} b, t = \partial^{-}_{i+1} a$.

Thus, in all cases, $c$ is $k$-thin, so the induction is complete.
The case $j=n-1$ of (*) completes the proof of Proposition 2.2.
\sqbox

\begin{Corollary} \label{C:2.7}
Let $\C$ be a cubical \w-category (or $m$-category) with
connections.
\begin{enumerate} [\rm (i)]
\item $n$-shells of the form $\bvarepsilon_i c$ or
$\mathbf{\Gamma}^{\alpha}_i c$ for $c \in \C_{n-1}$ are
commutative. \item Composites of commutative shells are
commutative.\sqbox
\end{enumerate}
\end{Corollary}

From Proposition 2.2 we easily deduce
\begin{theorem} \label{T:2.8}
Let $\C$ be a cubical \w-category (or $m$-category) with
connections. An element $c$ in $\C_n$ is thin if and only if it is
a composite of elements of the form $\varepsilon_i a$ or
$\Gamma^{\alpha}_i a \;(a \in \C_{n-1})$.
\end{theorem}

\Proof Proposition 2.2 shows that any such composite is thin. For
the converse, suppose that $c \in \C_n$ is thin. Then $\Psi c =
\psi_1 \psi_2 \ldots \psi_{n-1} c = \varepsilon_1 z$ for some $z
\in \C_{n-1}$. Now we saw in the proof of Lemma 1.5 that any
element $x \in \C_n$ can be written as a composite of $\psi_j x$
and elements of type $\varepsilon_i a, \Gamma^{\alpha}_i a$,
namely
$$
x = \quad \vcenter{ \xymatrix@M=0pt @=0.7pc{ \ar @{-} [dddddd] \ar
@{-} [rrrrrr] &&\ar @{-}
[dd] && \ar @{-} [dd] && \ar @{-} [dddddd] \\
& \sq  & & \vv  & & \tl & \\
\ar @{-} [rrrrrr] && \ar @{.} [dd] && \ar @{.} [dd] && \\
& \tl& & x & & \br &  \\
\ar @{-} [rrrrrr] && \ar @{-} [dd] && \ar @{-} [dd] && \\
& \br& & \vv  & & \sq  &  \\
\ar @{-} [rrrrrr]  &&&&&&}}  \quad  \directs{j+1}{j}
$$
(here the dotted segments indicate that $\psi_j x$ is first
partitioned as shown and the $3 \times 3$ array is then completed
so as to be composable.)

 By iteration, $x$ can be written as a
composite of $\varepsilon_1 z$ and elements of type $\varepsilon_i
a, \Gamma^{\alpha}_i a$. \sqbox

\begin{Corollary} \label{C:2.9}
An $n$-shell is commutative if and only if it can be written as a
composite of shells of type $\bvarepsilon_i a,
\mathbf{\Gamma}^{\alpha}_i a$, where $a \in \C_{n-1}$. \sqbox
\end{Corollary}
\noindent \textbf{Remark 1.} It is not clear from the proof of 2.8
whether a thin element can always be written as a composite of an
{\it array} of elements of type $\varepsilon_i a,
\Gamma^{\alpha}_i a$.\\

\noindent \textbf{Remark 2.} The particular folding map $\Psi$
used to define thin elements depends on a number of choices and
conventions. Theorem 2.8 shows that the notion of thinness is
intrinsic and does not depend on these choices. Thus, for example,
one might use $\Psi' = \psi_{n-1} \psi_{n-2} \ldots \psi_1$
instead of $\Psi$ but, by symmetry and Theorem 2.8, this would
give the same concept. Similarly, the more complicated full
folding operation $\Phi_n$ used in \cite{A-B-S02} gives the same
concept of thinness (see Section 9 of that paper, especially
Proposition 9.2 and Theorem 9.3). It is particularly reassuring
that the concept of commutative shell is independent of the choice
of foldings.

\section{Thin structures and connections}
We now consider cubical \w-categories (or $m$-categories) without
the assumption of the extra structure of connections. Of course
elements $\varepsilon_i a$ and shells $\bvarepsilon_i a $ exist
for such \w-categories so, in view of Theorem 2.8, it is not
surprising that there is a close relationship between existence of
thin elements and the existence of connections. An equivalence
between them in the 2-dimensional case was proved in \cite{B-M99}.
We extend this result to all dimensions.

Let $\C$ be a cubical \w-category (or $m$-category). Suppose that
$\C$ has connections $\; \Gamma^{-}_i, \Gamma^{+}_i: \C_{k-1} \to
\C_k$, defined for all $k=1,2,\ldots,n-1$, satisfying the usual
laws, up to that dimension (see \cite{A-B-S02}.) We aim to
characterize possible definitions of thin elements in $\C_n$
without first introducing more connections there.

As mentioned in Section 1, the $n$-category $(\C_0, \C_1, \ldots,
\C_{n-1}, \square \, \C_{n-1})$ does have connections in dimension
$n$ as well as those in lower dimension, so we can define folding
operations $\psi_1, \psi_2, \ldots, \psi_{n-1}, \Psi: \square \,
\C_{n-1} \to \square \, \C_{n-1}$. As a result, the idea of a
commutative $n$-shell is available, and we denote by $\nsq \,
\C_{n-1}$ the set of commutative $n$-shells in $\square \,
\C_{n-1}$. Clearly, by Corollary 2.7, $(\C_0, \C_1, \ldots,
\C_{n-1}, \nsq \, \C_{n-1})$ is a sub-(cubical $n$-category with
connections) of $(\C_0, \C_1, \ldots, \C_{n-1}, \square \,
\C_{n-1}).$\\

\noindent\textbf{Definition } A \emph{thin structure} on $\C_n$ is
a morphism $$\theta: (\C_0, \C_1, \ldots, \C_{n-1}, \nsq \,
\C_{n-1}) \to (\C_0, \C_1, \ldots, \C_{n-1}, \C_n)$$ of cubical
$n$-categories which is the identity on $\C_0, \C_1,
\ldots,\C_{n-1}$. Such a thin structure defines ``thin" elements
in $\C_n$, namely elements of the form $\theta(\mathbf{s})$ for a
commutative $n$-shell $\s$. Note that $\theta$ is necessarily
injective on $\nsq \, \C_{n-1}$ (because it preserves faces) and
the image of $\theta$ must be a sub-$n$-category of $\C$.
Consequently, every commutative $n$-shell in $\C$ has a unique
``thin" filler in $\C_n$, and the composites of ``thin" elements
are ``thin". Furthermore, we may now define $\Gamma^{\alpha}_i:
\C_{n-1} \to \C_n$ by $\Gamma^{\alpha}_i a = \theta
\bGamma^{\alpha}_i a$. Because $\theta$ preserves the lower
dimensional $\Gamma^{\alpha}_i$,  $ \o_i, \;
\partial^{\alpha}_i$ and $\varepsilon_i$, these newly defined
$\Gamma^{\alpha}_i$ satisfy the required laws making $(\C_0, \C_1,
\ldots, \C_{n-1}, \C_n)$ a cubical $n$-category with connections.
The thin elements of $\C_n$ defined using these
$\Gamma^{\alpha}_i$ are precisely the same as the ``thin" elements
defined by $\theta$, because of Proposition 2.1.

Conversely, if we are given connections $\Gamma^{\alpha}_i:
\C_{n-1} \to \C_n$ making $(\C_0, \C_1, \ldots, \C_n)$ a cubical
$n$-category with connections, there is a unique thin structure
$\theta$ with $\theta \bGamma^{\alpha}_i a = \Gamma^{\alpha}_i a$
for all $a \in \C_{n-1}, \alpha \in \{+,-\}, i \in
\{1,2,\ldots,n-1\}$. This is because the morphism $\theta$ must
map $\bvarepsilon_i a$ to $\varepsilon_i a$ and therefore, when it
is defined on the $\bGamma^{\alpha}_i a$, it is uniquely
determined on all commutative shells, by Corollary 2.9. That such
a $\theta$ exists is easily deduced from Proposition 2.1, and the
thin elements defined by the two methods again coincide. Hence
\begin{theorem} \label{T:3.1}
Let $\C = (\C_0, \C_1, \ldots, \C_n)$ be a cubical $n$-category
and suppose that $$(\C_0, \C_1, \ldots, \C_{n-1})$$ has the
structure of cubical $(n-1)$-category with connections. Then there
is a natural bijection between thin structures $\theta: \nsq
\,\C_{n-1} \to \C_n$ and sets of connections $\Gamma^{+}_i,
\Gamma^{-}_i: \C_{n-1} \to \C_n$ making $(\C_0, \C_1, \ldots,
\C_n)$ a cubical $n$-category with connections. The thin elements
defined by the connections coincide in all cases with the thin
elements defined by the corresponding $\theta$.
\end{theorem}

\noindent \textbf{Remark } It would be useful to have a simple
description of what is meant by a cubical $T$-complex, that is  a
weak thin structure (in all dimensions) on a cubical
\emph{complex}. The aim would be to impose axioms on the set of
``thin" elements in each dimension which would be equivalent to
the existence of a cubical \w-category structure with connections.
A simple description does exist in the groupoid  case (see
\cite{D77,B-H81a}), but seems to be more difficult in the category
case.
\medskip

\noindent {\bf Acknowledgement} It is a pleasure to record my
indebtedness to  Ronnie Brown for many discussions on the content
of this paper,  for support of its typing through his Leverhulme
Emeritus Fellowship, and for care to the production of the
diagrams in xypic. Our long and fruitful collaboration in the area
of higher dimensional algebra and homotopy has been a major source
of inspiration to me over the years.

\end{document}